\documentclass{amsart}
\usepackage{amstext}
\usepackage{amsthm}
\usepackage{amssymb}
\usepackage{esint}

\makeatletter
\numberwithin{equation}{section}
\numberwithin{figure}{section}
\theoremstyle{plain}
\newtheorem{thm}{Theorem}

\allowdisplaybreaks
\usepackage[sort&compress,square,numbers]{natbib}

\makeatother

\begin{document}

\title{Differential equations associated with Legendre polynomials}

\author{Taekyun Kim}
\address{Department of Mathematics, Kwangwoon University, Seoul 139-701, Republic
of Korea}
\email{tkkim@kw.ac.kr}

\author{Dae San Kim}
\address{Department of Mathematics, Sogang University, Seoul 121-742, Republic
of Korea}
\email{dskim@sogang.ac.kr}

\begin{abstract}
In this paper, we study non-linear differential equations associated with Legendre
polynomials and their applications. From our study of non-linear differential
equations, we derive some new and explicit identities for Legendre polynomials. 
\end{abstract}



\maketitle
\global\long\def\relphantom#1{\mathrel{\phantom{{#1}}}}

\section{Introduction}

The Legendre differential equation is given by 
\begin{equation}
\left(1-x^{2}\right)\frac{d^{2}y}{dx^{2}}-2x\frac{dy}{dx}+n\left(n+1\right)y=0,\quad\left(\text{see \cite{key-2,key-4,key-34}}\right).\label{eq:1}
\end{equation}

The equation (\ref{eq:1}) is equivalent to 
\begin{equation}
\frac{d}{dx}\left[\left(1-x^{2}\right)\frac{dy}{dx}\right]+n\left(n+1\right)y=0,\quad\left(\text{see \cite{key-7,key-21}}\right).\label{eq:2}
\end{equation}

The Legendre polynomials (or Legendre functions) are defined as the
solutions of Legendre differential equation. 

The generating function of Legendre polynomials $p_{n}\left(x\right)$
is given by 
\begin{equation}
\frac{1}{\sqrt{1-2xt+t^{2}}}=\sum_{n=0}^{\infty}p_{n}\left(x\right)t^{n},\quad\left(\text{see \cite{key-1,key-2,key-34}}\right).\label{eq:3}
\end{equation}

In physics, the generating function of Legendre polynomials is the
basis for multiple expansion. 

It is known that the Rodrigues' formula of Legendre polynomials is
given by 
\begin{equation}
p_{n}\left(x\right)=\frac{1}{2^{n}n!}\frac{d^{n}}{dx^{n}}\left[\left(x^{2}-1\right)^{n}\right],\quad\left(\text{see \cite{key-34}}\right).\label{eq:4}
\end{equation}

Thus, from (\ref{eq:4}), we note that $p_{n}\left(x\right)$ are
polynomials of degree $n$. 

The Legendre polynomials
$p_{n}\left(x\right)$ can also be represented by the contour integral
as 
\begin{equation}
p_{n}\left(x\right)=\frac{1}{2\pi i}\oint\left(1-2tx+t^{2}\right)^{-\frac{1}{2}}t^{-n-1}dt,\quad\left(\text{see \cite{key-2,key-21}}\right),\label{eq:5}
\end{equation}
where the contour encloses the origin and is traversed in a counterclockwise
direction. 

The first few Legendre polynomials are 
\begin{align*}
p_{0}\left(x\right) & =1,\quad p_{1}\left(x\right)=x,\quad p_{2}\left(x\right)=\frac{1}{2}\left(3x^{2}-1\right),\\
p_{3}\left(x\right) & =\frac{1}{2}\left(5x^{3}-3x\right),\\
p_{4}\left(x\right) & =\frac{1}{8}\left(35x^{4}-30x^{2}+3\right),\quad p_{5}\left(x\right)=\frac{1}{8}\left(63x^{5}-70x^{3}+15x\right),\\
p_{6}\left(x\right) & =\frac{1}{16}\left(231x^{6}-315x^{4}+105x^{2}-5\right),\quad\cdots.
\end{align*}

As is well known, the double factorial of a positive integer $n$
is a generalization of the usual factorial $n!$ defined by Arfken
and given by
\begin{equation}
n!!=\begin{cases}
n\left(n-2\right)\cdots5\cdot3\cdot1 & \text{ if }n>0\text{ odd},\\
n\left(n-2\right)\cdots6\cdot4\cdot2 & \text{if }n>0\text{ even},\\
1 & \text{if }n=-1,0,\quad\left(\text{see \cite{key-21}}\right)
\end{cases}.\label{eq:5-1}
\end{equation}

As was shown by  Arfken, the double factorial can be expressed in terms of the gamma function

\begin{equation}
\Gamma\left(n+\frac{1}{2}\right)=\frac{\left(2n-1\right)!!}{2^{n}}\sqrt{\pi},\quad\left(\text{see \cite{key-8,key-19,key-21}}\right).\label{eq:6}
\end{equation}

Thus, we note that the double factorial can also be extended to negative
odd integers using the definition: 
\begin{equation}
\left(-2n-1\right)!!=\frac{\left(-1\right)^{n}}{\left(2n-1\right)!!}=\frac{\left(-1\right)^{n}2^{n}n!}{\left(2n\right)!}.\label{eq:7}
\end{equation}

Now, we define the higher-order Legendre polynomials as follows: 
\begin{equation}
\left(\frac{1}{\sqrt{1-2tx+t^{2}}}\right)^{\alpha}=\sum_{n=0}^{\infty}p_{n}^{\left(\alpha\right)}\left(x\right)t^{n}.\label{eq:8}
\end{equation}

Some of the explicit formulas for $p_{n}\left(x\right)$ are 
\begin{align*}
p_{n}\left(x\right) & =\frac{1}{2^{n}}\sum_{k=0}^{n}\binom{n}{k}^{2}\left(x-1\right)^{n-k}\left(x+1\right)^{k}\\
 & =\sum_{k=0}^{n}\binom{n}{k}\binom{-n-1}{k}\left(\frac{1-x}{2}\right)^{k}\\
 & =2^{n}\sum_{k=0}^{n}\binom{n}{k}\binom{\frac{n+k-1}{2}}{n},\quad\left(\text{see \cite{key-2,key-34}}\right).
\end{align*}

Legendre polynomials occur in the solution of Laplacian equation of
the static potential $\nabla^{2}\phi\left(x\right)=0$, in a charge-free
region of space, using the method of separation of variables, where
the boundary condition has axial symmetry (no dependence on an azimuthal
angle).

Where $\hat{z}$ is the axis of symmetry and $\theta$ is the angle
between the position of the observer and the $\hat{z}$ axis (the
zenith angle), the solution for the potential will be 
\[
\phi\left(r,\theta\right)=\sum_{l=0}^{\infty}\left[A_{l}r^{l}+B_{l}r^{-\left(l+1\right)}\right]p_{l}\left(\cos\theta\right).
\]

Here $A_{l}$ and $B_{l}$ are to be determined according to the boundary
condition of each problem (see \cite{key-21}). 

Recently, several authors have studied some interesting extensions
and modifications of Legendre polynomials along with related combinatorial,
probabilistic, physics, and physical applications (see \cite{key-1,key-2,key-3,key-4,key-5,key-6,key-7,key-8,key-9,key-10,key-11,key-12,key-13,key-14,key-15,key-16,key-17,key-18,key-19,key-20,key-21,key-22,key-23,key-24,key-25,key-26,key-27,key-28,key-29,key-30,key-31,key-32,key-33}). 

Kim in \cite{key-25,key-26}, and Kim-Kim in \cite{key-24} considered some non-linear differential
equations arising from special numbers and polynomials and derived some new and interesting combinatorial identities. 

In this paper, we consider some differential equations arising from
the generating function of Legendre polynomials and give some new
and explicit identities on the Legendre polynomials which are derived
from the solutions of our differential equations.

\section{Differential equations arising from Legendre polynomials}

Let 
\begin{equation}
F=F\left(t,x\right)=\frac{1}{\sqrt{1-2tx+t^{2}}}.\label{eq:10}
\end{equation}

Then 
\begin{align}
F^{\left(1\right)} & =\frac{d}{dt}F\left(t,x\right)=\left(-\frac{1}{2}\right)\left(1-2tx+t^{2}\right)^{-\frac{3}{2}}\left(-2x+2t\right)\label{eq:11}\\
 & =\left(x-t\right)\left(1-2tx+t^{2}\right)^{-\frac{3}{2}}\nonumber \\
 & =\left(x-t\right)F^{3}.\nonumber 
\end{align}

From (\ref{eq:11}), we have 
\begin{equation}
F^{3}=\frac{1}{x-t}F^{\left(1\right)}.\label{eq:12}
\end{equation}

For $N\in\mathbb{N}$, let 
\begin{equation}
F^{\left(N\right)}=\left(\frac{d}{dt}\right)^{N}F\left(t,x\right),\quad F^{N}=\underset{N-\text{times}}{\underbrace{F\times\cdots\times F}}.\label{eq:13}
\end{equation}

From (\ref{eq:12}), we can derive the following equations: 
\begin{equation}
3F^{2}F^{\left(1\right)}=\frac{\left(-1\right)^{2}}{\left(x-t\right)^{2}}F^{\left(1\right)}+\frac{1}{x-t}F^{\left(2\right)}.\label{eq:14}
\end{equation}

By (\ref{eq:11}) and (\ref{eq:14}), we get 
\begin{equation}
3F^{5}=\frac{1}{\left(x-t\right)^{3}}F^{\left(1\right)}+\frac{1}{\left(x-t\right)^{2}}F^{\left(2\right)}.\label{eq:15}
\end{equation}

From (\ref{eq:15}), we have 
\begin{align}
3\cdot5F^{4}F^{\left(1\right)} & =\frac{3}{\left(x-t\right)^{4}}F^{\left(1\right)}+\frac{1}{\left(x-t\right)^{3}}F^{\left(2\right)}+\frac{2}{\left(x-t\right)^{3}}F^{\left(2\right)}+\frac{1}{\left(x-t\right)^{2}}F^{\left(3\right)}\label{eq:16}\\
 & =\frac{3}{\left(x-t\right)^{4}}F^{\left(1\right)}+\frac{3}{\left(x-t\right)^{3}}F^{\left(2\right)}+\frac{1}{\left(x-t\right)^{2}}F^{\left(3\right)}.\nonumber 
\end{align}

From (\ref{eq:11}) and (\ref{eq:16}), we get 
\begin{equation}
3\cdot5F^{7}=\frac{3}{\left(x-t\right)^{5}}F^{\left(1\right)}+\frac{3}{\left(x-t\right)^{4}}F^{\left(2\right)}+\frac{1}{\left(x-t\right)^{3}}F^{\left(3\right)}.\label{eq:17}
\end{equation}

From (\ref{eq:17}), we note that 
\begin{align}
3\cdot5\cdot7F^{6}F^{\left(1\right)} & =\frac{5\cdot3}{\left(x-t\right)^{6}}F^{\left(1\right)}+\frac{3}{\left(x-t\right)^{5}}F^{\left(2\right)}+\frac{4\cdot3}{\left(x-t\right)^{5}}F^{\left(2\right)}+\frac{3}{\left(x-t\right)^{4}}F^{\left(3\right)}\label{eq:18}\\
 & \relphantom =+\frac{3}{\left(x-t\right)^{4}}F^{\left(3\right)}+\frac{1}{\left(x-t\right)^{3}}F^{\left(4\right)}.\nonumber 
\end{align}

Thus, by (\ref{eq:11}) and (\ref{eq:18}), we get 
\begin{equation}
3\cdot5\cdot7F^{9}=\frac{5\cdot3}{\left(x-t\right)^{7}}F^{\left(1\right)}+\frac{5\cdot3}{\left(x-t\right)^{6}}F^{\left(3\right)}+\frac{3\cdot2}{\left(x-t\right)^{5}}F^{\left(3\right)}+\frac{1}{\left(x-t\right)^{4}}F^{\left(4\right)}.\label{eq:19}
\end{equation}

Continuing this process, we set
\begin{equation}
\left(2N-1\right)!!F^{2N+1}=\sum_{i=1}^{N}a_{i}\left(N\right)\frac{F^{\left(i\right)}}{\left(x-t\right)^{2N-i}}.\label{eq:20}
\end{equation}

Thus, by (\ref{eq:20}), we get 
\begin{align}
 & \left(2N-1\right)!!\left(2N+1\right)F^{2N}F^{\left(1\right)}\label{eq:21}\\
 & =\sum_{i=1}^{N}a_{i}\left(N\right)\frac{d}{dt}\left(\frac{F^{\left(i\right)}}{\left(x-t\right)^{2N-i}}\right)\nonumber \\
 & =\sum_{i=1}^{N}a_{i}\left(N\right)\left\{ \frac{2N-i}{\left(x-t\right)^{2N-i+1}}F^{\left(i\right)}+\frac{1}{\left(x-t\right)^{2N-i}}F^{\left(i+1\right)}\right\} \nonumber \\
 & =\sum_{i=1}^{N}a_{i}\left(N\right)\frac{2N-i}{\left(x-t\right)^{2N-i+1}}F^{\left(i\right)}+\sum_{i=1}^{N}a_{i}\left(N\right)\frac{F^{\left(i+1\right)}}{\left(x-t\right)^{2N-i}}.\nonumber 
\end{align}

From (\ref{eq:11}) and (\ref{eq:21}), we have 
\begin{align}
 & \left(2N+1\right)!!F^{2N+3}\label{eq:22}\\
 & =\sum_{i=1}^{N}a_{i}\left(N\right)\frac{2N-i}{\left(x-t\right)^{2\left(N+1\right)-i}}F^{\left(i\right)}+\sum_{i=1}^{N}a_{i}\left(N\right)\frac{F^{\left(i+1\right)}}{\left(x-t\right)^{2N-i+1}}\nonumber \\
 & =\sum_{i=1}^{N}a_{i}\left(N\right)\frac{2N-i}{\left(x-t\right)^{2\left(N+1\right)-i}}F^{\left(i\right)}+\sum_{i=1}^{N}a_{i}\left(N\right)\frac{F^{\left(i+1\right)}}{\left(x-t\right)^{2N-i+1}}\nonumber \\
 & =\left(2N-1\right)a_{1}\left(N\right)\frac{F^{\left(1\right)}}{\left(x-t\right)^{2N+1}}+a_{N}\left(N\right)\frac{F^{\left(N+1\right)}}{\left(x-t\right)^{N+1}}\nonumber \\
 & \relphantom =+\sum_{i=2}^{N}\left\{ \left(2N-i\right)a_{i}\left(N\right)+a_{i-1}\left(N\right)\right\} \frac{F^{\left(i\right)}}{\left(x-t\right)^{2\left(N+1\right)-i}}.\nonumber 
\end{align}

By replacing $N$ by $N+1$ in (\ref{eq:20}), we get 
\begin{equation}
\left(2N+1\right)!!F^{2N+3}=\sum_{i=1}^{N+1}a_{i}\left(N+1\right)\frac{F^{\left(i\right)}}{\left(x-t\right)^{2\left(N+1\right)-i}}.\label{eq:23}
\end{equation}

By comparing the coefficients on the both sides of (\ref{eq:22})
and (\ref{eq:23}), we have 
\begin{align}
a_{1}\left(N+1\right) & =\left(2N-1\right)a_{1}\left(N\right),\label{eq:24}\\
a_{N+1}\left(N+1\right) & =a_{N}\left(N\right),\label{eq:25}
\end{align}
and 
\begin{equation}
a_{i}\left(N+1\right)=\left(2N-i\right)a_{i}\left(N\right)+a_{i-1}\left(N\right),\quad\left(2\le i\le N\right).\label{eq:26}
\end{equation}

From (\ref{eq:12}) and (\ref{eq:20}), we note that 
\begin{equation}
\frac{1}{x-t}F^{\left(1\right)}=F^{3}=a_{1}\left(1\right)\frac{1}{x-t}F^{\left(1\right)}.\label{eq:27}
\end{equation}

Thus, by comparing the coefficients on both sides of (\ref{eq:27}),
we get $a_{1}\left(1\right)=1$. 

From (\ref{eq:24}) and (\ref{eq:25}), we can derive the following
equations: 
\begin{align}
a_{1}\left(N+1\right) & =\left(2N-1\right)a_{1}\left(N\right)\label{eq:28}\\
 & =\left(2N-1\right)\left(2N-3\right)a_{1}\left(N-1\right)\nonumber \\
 & \vdots\nonumber \\
 & =\left(2N-1\right)\left(2N-3\right)\cdots3\cdot1a_{1}\left(1\right)\nonumber \\
 & =\left(2N-1\right)!!,\nonumber 
\end{align}
and 
\begin{equation}
a_{N+1}\left(N+1\right)=a_{N}\left(N\right)=a_{N-1}\left(N-1\right)=\cdots=a_{1}\left(1\right)=1.\label{eq:29}
\end{equation}

By (\ref{eq:26}), for $2\le i \le N$, we get 
\begin{align}
 & a_{i}\left(N+1\right)\label{eq:30}\\
 & =\left(2N-i\right)a_{i}\left(N\right)+a_{i-1}\left(N\right)\nonumber \\
 & =\left(2N-i\right)\left\{ \left(2\left(N-1\right)-i\right)a_{i}\left(N-1\right)+a_{i-1}\left(N-1\right)\right\} +a_{i-1}\left(N\right)\nonumber \\
 & =\left(2N-i\right)\left(2N-2-i\right)a_{i}\left(N-1\right)+\left(2N-i\right)a_{i-1}\left(N-1\right)+a_{i-1}\left(N\right)\nonumber \\
 & =\left(2N-i\right)\left(2N-2-i\right)\left(2N-4-i\right)a_{i}\left(N-2\right)+\left(2N-i\right)\left(2N-2-i\right)a_{i-1}\left(N-2\right)\nonumber \\
 & \relphantom =+\left(2N-i\right)a_{i-1}\left(N-1\right)+a_{i-1}\left(N\right)\nonumber \\
 & \vdots\nonumber \\
 & =\left(\prod_{l=0}^{N-i}\left(2N-2l-i\right)\right)a_{i}\left(i\right)+\sum_{l=0}^{N-i}\left\langle 2N-i\right\rangle _{l}a_{i-1}\left(N-l\right)\nonumber \\
 & =\prod_{l=0}^{N-i}\left(2N-2l-i\right)+\sum_{l=0}^{N-i}\left\langle 2N-i\right\rangle _{l}a_{i-1}\left(N-l\right),\nonumber 
\end{align}
where $\left\langle 2N+\alpha\right\rangle _{k}=\left(2N+\alpha\right)\left(2\left(N-1\right)+\alpha\right)\cdots\left(2\left(N-k+1\right)+\alpha\right),$
and $\left\langle 2N+\alpha\right\rangle _{0}=1$. 

From (\ref{eq:30}), we have 
\begin{align}
 & a_{i-1}\left(N-l_{1}\right)\label{eq:31}\\
 & =\prod_{l_{2}=0}^{N-l_{1}-i}\left(2N-2l_{1}-2l_{2}-i-1\right)\nonumber \\
 & \relphantom =+\sum_{l_{2}=0}^{N-l_{1}-i}\left\langle 2N-2l_{1}-i-1\right\rangle _{l_{2}}a_{i-2}\left(N-l_{1}-l_{2}-1\right).\nonumber 
\end{align}

By (\ref{eq:30}) and (\ref{eq:31}), we get 
\begin{align}
 & a_{i}\left(N+1\right)\label{eq:32}\\
 & =\prod_{l_{1}=0}^{N-i}\left(2N-2l_{1}-i\right)+\sum_{l_{1}=0}^{N-i}\left\langle 2N-i\right\rangle _{l_{1}}a_{i-1}\left(N-l_{1}\right)\nonumber \\
 & =\prod_{l_{1}=0}^{N-i}\left(2N-2l_{1}-i\right)+\sum_{l_{1}=0}^{N-i}\left(\prod_{l_{2}=0}^{N-l_{1}-i}\left(2N-2l_{1}-2l_{2}-i-1\right)\right)\left\langle 2N-i\right\rangle _{l_{1}}\nonumber \\
 & \relphantom =+\sum_{l_{1}=0}^{N-i}\sum_{l_{2}=0}^{N-l_{1}-i}\left\langle 2N-i\right\rangle _{l_{1}}\left\langle 2N-2l_{1}-i-1\right\rangle _{l_{2}}a_{i-2}\left(N-l_{1}-l_{2}-1\right).\nonumber 
\end{align}

Now, we observe that 
\begin{align}
 & a_{i-2}\left(N-l_{1}-l_{2}-1\right)\label{eq:33}\\
 & =a_{i-2}\left(N-l_{1}-l_{2}-2+1\right)\nonumber \\
 & =\prod_{l_{3}=0}^{N-l_{1}-l_{2}-i}\left(2N-2l_{1}-2l_{2}-2l_{3}-i+2\right)\nonumber \\
 & \relphantom =+\sum_{l_{3}=0}^{N-l_{1}-l_{2}-i}\left\langle 2N-2l_{1}-2l_{2}-i-2\right\rangle _{l_{3}}a_{i-3}\left(N-l_{1}-l_{2}-l_{3}-2\right).\nonumber 
\end{align}

From (\ref{eq:32}) and (\ref{eq:33}), and continuing this process,  we obtain
\begin{align*}
 & \relphantom =a_{i}\left(N+1\right)\\
 & =\prod_{l_{1}=0}^{N-i}\left(2N-2l_{1}-i\right)\\
 & \relphantom =+\sum_{l_{1}=0}^{N-i}\left\langle 2N-i\right\rangle _{l_{1}}\left(\prod_{l_{2}=0}^{N-l_{1}-i}\left(2N-2l_{1}-2l_{2}-i-1\right)\right)\\
 & \relphantom =+\sum_{l_{1}=0}^{N-i}\sum_{l_{2}=0}^{N-l_{1}-i}\left\langle 2N-i\right\rangle _{l_{1}}\left\langle 2N-2l_{1}-i-1\right\rangle _{l_{2}}\\
 & \relphantom =\times\left(\prod_{l_{3}=0}^{N-l_{1}-l_{2}-i}\left(2N-2l_{1}-2l_{2}-2l_{3}-i-2\right)\right)+\cdots\\
 & \relphantom =+\sum_{l_{1}=0}^{N-i}\sum_{l_{2}=0}^{N-l_{1}-i}\cdots\sum_{l_{i-2}=0}^{N-l_{1}-\cdots-l_{i-3}-i}\left\langle 2N-i\right\rangle _{l_{1}}\\
 & \relphantom =\times\left\langle 2N-2l_{1}-i-1\right\rangle _{l_{2}}\cdots\left\langle 2N-2l_{1}-\cdots-2l_{i-3}-2i+3\right\rangle _{l_{i-2}}\\
 & \relphantom =\times\left(\prod_{l_{i-1}=0}^{N-l_{1}-\cdots-l_{i-2}-i}\left(2N-2l_{1}-2l_{2}-\cdots-2l_{i-1}-2i+2\right)\right)\\
 & \relphantom =+\sum_{l_{1}=0}^{N-i}\sum_{l_{2}=0}^{N-l_{1}-i}\cdots\sum_{l_{i-1}=0}^{N-l_{1}-l_{2}-\cdots-l_{i-2}-i}\left\langle 2N-i\right\rangle _{l_{1}}\\
 & \relphantom =\times\left\langle 2N-2l_{1}-i-1\right\rangle _{l_{2}}\cdots\left\langle 2N-2l_{1}-2l_{2}-\cdots-2l_{i-2}-2i+2\right\rangle _{l_{i-1}}\\
 & \relphantom =\times a_{1}\left(N-l_{1}-l_{2}-\cdots-l_{i-1}-i+2\right)\nonumber \\
 & =\prod_{l_{1}=0}^{N-i}\left(2N-2l_{1}-i\right)\\
 & \relphantom =+\sum_{l_{1}=0}^{N-i}\left\langle 2N-i\right\rangle _{l_{1}}\left(\prod_{l_{2}=0}^{N-l_{1}-i}\left(2N-2l_{1}-2l_{2}-i-1\right)\right)\\
 & \relphantom =+\sum_{l_{1}=0}^{N-i}\sum_{l_{2}=0}^{N-l_{1}-i}\left\langle 2N-i\right\rangle _{l_{1}}\left\langle 2N-2l_{1}-i-1\right\rangle _{l_{2}}\\
 & \relphantom =\times\left(\prod_{l_{3}=0}^{N-l_{1}-l_{2}-i}\left(2N-2l_{1}-2l_{2}-2l_{3}-i-2\right)\right)+\cdots\\
 & \relphantom =+\sum_{l_{1}=0}^{N-i}\sum_{l_{2}=0}^{N-l_{1}-i}\cdots\sum_{l_{i-2}=0}^{N-l_{1}-\cdots-l_{i-3}-i}\left\langle 2N-i\right\rangle _{l_{1}}\\
 & \relphantom =\times\left\langle 2N-2l_{1}-i-1\right\rangle _{l_{2}}\cdots\left\langle 2N-2l_{1}-\cdots-2l_{i-3}-2i+3\right\rangle _{l_{i-2}}\\
 & \relphantom =\times\left(\prod_{l_{i-1}=0}^{N-l_{1}-\cdots-l_{i-2}-i}\left(2N-2l_{1}-2l_{2}-\cdots-2l_{i-1}-2i+2\right)\right)\\
 & \relphantom =+\sum_{l_{1}=0}^{N-i}\sum_{l_{2}=0}^{N-l_{1}-i}\cdots\sum_{l_{i-1}=0}^{N-l_{1}-l_{2}-\cdots-l_{i-2}-i}\left\langle 2N-i\right\rangle _{l_{1}}\\
 & \relphantom =\times\left\langle 2N-2l_{1}-i-1\right\rangle _{l_{2}}\cdots\left\langle 2N-2l_{1}-2l_{2}-\cdots-2l_{i-2}-2i+2\right\rangle _{l_{i-1}}\\
 & \relphantom =\times\left(2\left(N-l_{1}-l_{2}-\cdots-l_{i-1}-i\right)+1\right)!!.
\end{align*}

Therefore, we obtain the following theorem.
\begin{thm}
\label{thm:1} The following non-linear differential equations
\[
\left(2N-1\right)!!F^{2N+1}=\sum_{i=1}^{N}a_{i}\left(N\right)\frac{F^{\left(i\right)}}{\left(x-t\right)^{2N-i}},\quad\left(N=1,2,\dots\right),
\]
has a solution 
\[
F=F\left(t,x\right)=\frac{1}{\sqrt{1-2tx+t^{2}}},
\]
where $a_{1}\left(N\right)=\left(2N-3\right)!!$, $a_{N}\left(N\right)=1$, and, for $2 \le i \le N-1$,
\begin{align*}
 & \relphantom =a_{i}\left(N\right)\\
 & =\prod_{l_{1}=0}^{N-1-i}\left(2N-2l_{1}-i-2\right)\\
 & \relphantom =+\sum_{l_{1}=0}^{N-i-1}\left\langle 2N-2-i\right\rangle _{l_{1}}\left(\prod_{l_{2}=0}^{N-l_{1}-i-1}\left(2N-2l_{1}-2l_{2}-i-3\right)\right)\\
 & \relphantom =+\sum_{l_{1}=0}^{N-i-1}\sum_{l_{2}=0}^{N-l_{1}-i-1}\left\langle 2N-i-2\right\rangle _{l_{1}}\left\langle 2N-2l_{1}-i-3\right\rangle _{l_{2}}\\
 & \relphantom =\times\left(\prod_{l_{3}=0}^{N-l_{1}-l_{2}-i-1}\left(2N-2l_{1}-2l_{2}-2l_{3}-i-4\right)\right)+\cdots\\
 & \relphantom =+\sum_{l_{1}=0}^{N-i-1}\sum_{l_{2}=0}^{N-l_{1}-i-1}\cdots\sum_{l_{i-2}=0}^{N-l_{1}-\cdots-l_{i-3}-i-1}\left\langle 2N-i-2\right\rangle _{l_{1}}\\
 & \relphantom =\times\left\langle 2N-2l_{1}-i-3\right\rangle _{l_{2}}\cdots\left\langle 2N-2l_{1}-\cdots-2l_{i-3}-2i+1\right\rangle _{l_{i-2}}\\
 & \relphantom =\times\left(\prod_{l_{i-1}=0}^{N-l_{1}-\cdots-l_{i-2}-i-1}\left(2N-2l_{1}-2l_{2}-\cdots-2l_{i-1}-2i\right)\right)\\
 & \relphantom =+\sum_{l_{1}=0}^{N-i-1}\sum_{l_{2}=0}^{N-l_{1}-i-1}\cdots\sum_{l_{i-1}=0}^{N-l_{1}-l_{2}-\cdots-l_{i-2}-i-1}\left\langle 2N-i-2\right\rangle _{l_{1}}\\
 & \relphantom =\times\left\langle 2N-2l_{1}-i-3\right\rangle _{l_{2}}\cdots\left\langle 2N-2l_{1}-2l_{2}-\cdots-2l_{i-2}-2i\right\rangle _{l_{i-1}}\\
 & \relphantom =\times\left(2\left(N-l_{1}-l_{2}-\cdots-l_{i-1}-i\right)-1\right)!!.
\end{align*}
\end{thm}
Recall that the generating function of Legendre polynomials is given
by 
\begin{equation}
F=F\left(t,x\right)=\frac{1}{\sqrt{1-2tx+t^{2}}}=\sum_{n=0}^{\infty}p_{n}\left(x\right)t^{n}.\label{eq:36}
\end{equation}

Thus, by (\ref{eq:36}), we get 
\begin{align}
F^{\left(i\right)} & =\left(\frac{d}{dt}\right)^{i}F\left(t,x\right)=\sum_{n=i}^{\infty}p_{n}\left(x\right)\left(n\right)_{i}t^{n-i}\label{eq:37}\\
 & =\sum_{n=0}^{\infty}p_{n+i}\left(x\right)\left(n+i\right)_{i}t^{n},\nonumber 
\end{align}
where $\left(x\right)_{n}=x\left(x-1\right)\cdots\left(x-n+1\right)$,
$\left(n\ge1\right)$, and $\left(x\right)_{0}=1$. 

From (\ref{eq:8}), we note that 
\begin{equation}
F^{2N+1}=\sum_{n=0}^{\infty}p_{n}^{\left(2N+1\right)}\left(x\right)t^{n}.\label{eq:38}
\end{equation}

By Theorem \ref{thm:1}, we get 
\begin{align}
 & \relphantom =F^{2N+1}\label{eq:39}\\
 & =\frac{1}{\left(2N-1\right)!!}\sum_{i=1}^{N}a_{i}\left(N\right)\frac{F^{\left(i\right)}}{\left(x-t\right)^{2N-i}}\nonumber \\
 & =\frac{1}{\left(2N-1\right)!!}\sum_{i=1}^{N}a_{i}\left(N\right)\left(\sum_{m=0}^{\infty}\binom{2N+m-i-1}{m}x^{-\left(2N+m-i\right)}t^{m}\right)\nonumber \\
 & \relphantom =\times\left(\sum_{l=0}^{\infty}p_{l+i}\left(x\right)\left(l+i\right)_{i}t^{l}\right)\nonumber \\
 & =\frac{1}{\left(2N-1\right)!!}\sum_{i=1}^{N}a_{i}\left(N\right)\nonumber \\
 & \relphantom =\times\sum_{n=0}^{\infty}\left\{ \sum_{m=0}^{n}\binom{2N+m-i-1}{m}x^{-\left(2N+m-i\right)}p_{n-m+i}\left(x\right)\left(n-m+i\right)_{i}\right\} t^{n}.\nonumber 
\end{align}

Therefore, by (\ref{eq:38}) and (\ref{eq:39}), we obtain the following
theorem.
\begin{thm}
\label{thm:2} For $n\ge0$, we have 
\begin{align*}
p_{n}^{\left(2N+1\right)}\left(x\right) & =\frac{1}{\left(2N-1\right)!!}\sum_{i=1}^{N}\sum_{m=0}^{n}a_{i}\left(N\right)\binom{2N+m-i-1}{m}\\
 & \relphantom =\times x^{-\left(2N+m-i\right)}p_{n-m+i}\left(x\right)\left(n-m+i\right)_{i},
\end{align*}
where $a_{i}\left(N\right)$ $\left(1\le i \le N\right)$ are as in Theorem 1.
\end{thm}
\bibliographystyle{amsplain}
\providecommand{\bysame}{\leavevmode\hbox to3em{\hrulefill}\thinspace}
\providecommand{\MR}{\relax\ifhmode\unskip\space\fi MR }
\providecommand{\MRhref}[2]{%
  \href{http://www.ams.org/mathscinet-getitem?mr=#1}{#2}
}
\providecommand{\href}[2]{#2}

\end{document}